\documentclass[a4paper,12pt]{elsarticle}
\usepackage{fullpage}
\usepackage[utf8]{inputenc}
\usepackage[T1]{fontenc}
\usepackage{amsfonts}
\usepackage{mathrsfs,amsmath}
\usepackage{amssymb}
\usepackage{amsthm}
\usepackage{subfig}
\usepackage{mwe} 
\usepackage{multirow}
\usepackage{setspace}
\usepackage{here}
\usepackage{commath}
\usepackage{graphicx}
\usepackage{algorithm, algorithmic}
\usepackage{color}
\usepackage{verbatim}
\usepackage{subfig}
\usepackage{bigfoot} 
\usepackage{filecontents}
\usepackage{dcolumn}
\newtheorem{Teo}{Theorem}

\newdefinition{rmk}{Remark}

\DeclareMathOperator*{\re}{Re}
\DeclareMathOperator*{\im}{Im}
\DeclareMathOperator*{\cond}{cond}
\DeclareMathOperator*{\rexii}{REXII}
\DeclareMathOperator*{\rexi}{REXI}
\DeclareMathOperator*{\diag}{diag}

\begin{document}

\title{An accurate and time-parallel rational exponential integrator for hyperbolic and oscillatory PDEs\tnoteref{t1}} 
\tnotetext[t1]{The numerical computations are performed using the GPU cluster GPU3 at the University of Innsbruck.}

\author[1]{Marco Caliari}\ead{marco.caliari@univr.it}

\author[2]{Lukas Einkemmer}\ead{lukas.einkemmer@uibk.ac.at}

\author[2]{Alexander Moriggl}\ead{alexander.moriggl@uibk.ac.at}

\author[2]{Alexander Ostermann}\ead{alexander.ostermann@uibk.ac.at}

\address[1]{Department of Computer Science, University of Verona, Italy}
\address[2]{Department of Mathematics, University of Innsbruck, Austria}

\date{\today}

\begin{abstract}
    Rational exponential integrators (REXI) are a class of numerical methods that are well suited for the time integration of linear partial differential equations with imaginary eigenvalues. Since these methods can be parallelized in time (in addition to the spatial parallelization that is commonly performed) they are well suited to exploit modern high performance computing systems. In this paper, we propose a novel REXI scheme that drastically improves accuracy and efficiency. The chosen approach will also allow us to easily determine how many terms are required in the approximation in order to obtain accurate results. We provide comparative numerical simulations for a shallow water equation that highlight the efficiency of our approach and demonstrate that REXI schemes can be efficiently implemented on graphic processing units.
\end{abstract}


\maketitle

\section{Introduction}

In this work, we are interested in simulating linear partial differential equations (PDEs) with purely imaginary eigenvalues of large modulus (i.e.~stiff problems). That is, we consider 
\[
\partial_t f = Af, \quad f(0) = f_0, \quad \sigma(A) \subset i\mathbb{R},
\]
where $\sigma(A)$ is the spectrum of $A$. Such problems arise, for example, in quantum dynamics (e.g. the Schr\"odinger equation) and wave propagation (e.g. the Helmholtz equation). In addition, solving these problems is an integral part of applying exponential integrators or splitting methods to a large number of nonlinear problems ranging from plasma physics to electrodynamics.

Due to the stiff nature of these equations, explicit time stepping methods are forced to take excessively small time steps in order to remain stable. Thus, implicit schemes (e.g.~\cite{stabilityMultistepCrankNicolson,implicitTimeStepping}), implicit-explicit IMEX schemes (e.g.~\cite{ExpVSimex}), exponential integrators (e.g.~\cite{hochbruck_ostermann_2010, auer2018magnus, expAtmosphere,   crouseilles2020exponential,  crouseilles2018exponential}), or splitting methods (e.g.~\cite{ bao2002time,caliari2017splitting}) are commonly used. These methods enjoy better stability properties and can thus, in principle, take large time steps. However, for highly oscillatory problems, the maximal time step size of implicit methods is still severely limited by the fact that such methods need to resolve the oscillations.

Approximately a decade ago, the hardware used to run such simulations has undergone a paradigm shift. Due to the fact that frequency scaling has essentially ended at that point, the main way to increase performance has been to add more parallelism. Desktop computers now routinely have 16 cores and large vectorization units. This trend is even more pronounced in high performance computing systems, where supercomputers with millions of threads are now in operation. In addition, graphic processing units (GPUs) have come to the forefront as they are able to outperform central processing units (CPUs) for many scientific computing tasks. This advantage is achieved by providing a massive parallel system. In fact, a sequential program on a GPU would be slower than on a CPU. Therefore, methods which are highly parallelizable and work well on these new computer architectures are needed to take advantage of their computational power. A significant body of research has been accumulated in recent years that considers numerical methods that are well suited for such systems (see, e.g., \cite{einkemmer2017evaluation,einkemmer2020semi,guo2016,muller2013matrix,murray2011gpu}). More specifically, in the context of exponential integrators we refer to \cite{einkemmer2013exponential,farquhar2016}.

Since time stepping methods are inherently sequential, they generally can not be parallelized in time (although some parallelism can be extracted by constructing methods with parallel stages; see, e.g., \cite{luan2016parallel, nievergelt1964parallel}). Even exponential integrators of high order, which are able to perform large time steps, are mainly implemented sequentially (in time). For example, polynomial approximations, as in \cite{MohyHigham2011,LejaRevisted, Leja}, require us to calculate a sequence of matrix-vector products which can not be done in parallel. A similar argument holds for Krylov approximations~\cite{Krylov}. Of course, these methods are parallelizable in space. But in some situations the scalability is limited and, at some point, increasing the number of computing cores does not further reduce the simulation time~\cite{SPHW2017}. Therefore, such an approach can not fully exploit modern computer hardware.

A novel idea to overcome this problem are so-called Rational Exponential Integrators (REXI) schemes, which were introduced in \cite{HBMW2015}. The basic idea of these methods is to approximate $e^{tA}$ by a linear combination of simple rational functions. The advantage of REXI methods is that the corresponding terms can be calculated independently of each other. Therefore, these methods are highly parallelizable in time. It is worth mentioning that REXI is markedly distinct from time parallelization schemes such as the parareal and similar methods (see, e.g., \cite{gander201550}). In the latter case, a coarse and a fine time integrator are combined to achieve parallelism within an iterative procedure, while in the former the action of the matrix exponential is directly approximated in a way that is amendable to parallelization.

In this paper our goal is twofold. First, we propose a modification to the original REXI scheme that drastically improves the accuracy and efficiency of the method (section \ref{sec:REXI}). The proposed method also allows us to easily determine how many terms the approximation requires in order to obtain accurate results. These theoretical considerations are then confirmed by numerical experiments in section \ref{sec:NumericalExamples}. Second, we demonstrate that these types of methods can be efficiently implemented on massively parallel computer architectures. Specifically, we demonstrate an implementation on modern GPUs that yields a drastic speedup compared to the corresponding CPU implementation for the shallow water equations (section \ref{sec:numericalExampleLSWE}).

\section{The original and improved REXI schemes}
\label{sec:REXI}

In this section we discuss the derivation of the REXI schemes. Moreover, we reveal some problems of the original scheme in the matrix case and show how they can be eliminated with our new formulation.

\subsection{The scalar case}
\label{sec:REXI_derivation_scalar}
In this section we will give a brief summary of how REXI approximates $e^{ix}$ with $ x\in\mathbb{R}$. For more details, see \cite{HBMW2015, SPHW2017}. Our notation is the same as in \cite{SPHW2017}. The three main steps are as follows:
\begin{enumerate}
\itemsep0em 
\item Approximate $e^{ix}$ by a sum of Gaussian functions.
\item Approximate each Gaussian function by a sum of rational functions.
\item Combine 1. and 2. to approximate $e^{ix}$ by a sum of rational functions.
\end{enumerate}

\subsubsection*{Step 1}

We start by writing $e^{ix}$ as a linear combination of Gaussian functions \cite{GausKernelApprox}
\begin{equation}
\label{eq:eixsuminfty}
e^{ix} + \epsilon = f(x) = \sum_{m = -\infty}^\infty b_m \psi_h(x + mh),
\end{equation} 
where
\begin{equation}
\label{eq:psi}
\psi_h(x) = \frac{1}{\sqrt{4\pi}}\exp\left(-\frac{x^2}{4h^2}\right)
\end{equation}
and $\epsilon = e^{h^2}\sum_{k \neq 0}e^{-4\pi^2k^2}$, see Appendix B. Since $\epsilon$, in general, is smaller than machine precision (see Appendix B or \cite{GausKernelApprox}), we neglect it in the following.

The parameter $h$ defines the numerical support of the Gaussian function. Our goal now is to determine the coefficients $b_m$. To that end equation~\eqref{eq:eixsuminfty} is transformed to Fourier space. Using the shift property of the Fourier transform, it follows that
\[
\frac{\hat{f}(\omega)}{\hat{\psi}_h(\omega)} = \sum_{m = -\infty}^{\infty} b_m e^{2\pi i mh\omega}, 
\]
and therefore we obtain 
\begin{equation}
\label{eq:bmdefinitionintegral}
b_m = h\int_{-\frac{1}{2h}}^{\frac{1}{2h}} e^{-2\pi i m h \omega} \frac{\hat{f}(\omega)}{\hat{\psi}_h(\omega)} \, d\omega.
\end{equation}
Since $f(x) = e^{ix}$, we have $\hat{f}(\omega) = \delta(\omega - \frac{1}{2\pi})$, where $\delta(x)$ is the Dirac distribution. Moreover 
\begin{equation}
\label{eq:psi_fourier}
\hat{\psi}_h(\omega) = he^{-4\pi^2h^2\omega^2}
\end{equation}
and therefore the sought after coefficients are given by
\begin{equation}
\label{eq:bm}
b_m = e^{-imh}e^{h^2}.
\end{equation}
Here, the first constraint on $h$ arises. If $h > \pi$ then the integral in \eqref{eq:bmdefinitionintegral} is~$0$. In practice, this parameter has to be even smaller to produce accurate approximations. Its value is discussed in Section~\ref{sec:numericalExamplesScalar}.

Finally, the sum in \eqref{eq:eixsuminfty} has to be truncated:
\begin{equation}
\label{eq:eixsumM}
e^{ix} \approx \sum_{m = -M}^M b_m \psi_h(x + mh).
\end{equation} 
Therefore, REXI depends on two parameters: $M$ and $h$. As mentioned earlier, the parameter $h < \pi$ defines the numerical support of \eqref{eq:psi} and the parameter $M$ controls, together with $h$, the interval on which the approximation is sought. It can be shown (\cite{HBMW2015} and Appendix B), that \eqref{eq:eixsumM} produces an accurate approximation if 
\begin{equation}
\label{eq:Mformula}
|x| \leq (M-11)h.
\end{equation}

\subsubsection*{Step 2}

The second step consists in approximating the Gaussian function $\psi_h(x)$ as a sum of rational functions:

\begin{equation}
\label{eq:ratapproxgaus}
\psi_h(x) \approx \re\left(\sum_{l = -L}^{L}\frac{a_l}{i\frac{x}{h} + \mu + il}\right) =: R\left(\frac{x}{h}\right),
\end{equation}
where $\mu \in \mathbb{R}$ and $a_l \in \mathbb{C}$ are coefficients. To determine these parameters, the authors in \cite{HBMW2015} first approximate $\hat{\psi}_1(\omega)$ by a linear combination of exponential functions, 

\[
\hat{\psi}_1(\omega) \approx \sum_{j=1}^{J} b_j e^{\theta_j \omega},
\]
using the Adamyan--Arov--Krein theory (see \cite{DAMLE2013} for more details). Moving back to physical space, they then obtain
\begin{equation}
\label{eq:origAAK}
\mathscr{F}^{-1}\left(\sum_{j=1}^{J} b_j e^{\theta_j \omega}\right) = -2\re\left(\sum_{j=1}^{J} \frac{b_j}{2\pi ix+\theta_j}\right) \approx \psi_1(x),
\end{equation}
where $\re(\theta_j) < 0$. In the original $\rexi$ scheme \cite{HBMW2015} only the coefficients $\theta_j$ are computed. Then, they set $\mu = \min_j \re(\theta_j / 2\pi)$ and then look for an approximation of $\psi_1(x)$ of the form~\eqref{eq:ratapproxgaus} with $h=1$,

\begin{equation}
\label{eq:ratapproxgaush1}
R(x) = \re\left(\sum_{l = -L}^L \frac{a_l}{ix + \mu + il}\right),
\end{equation}
where $L = 11$ is chosen. This slight modification of \eqref{eq:origAAK} gives the advantage that all the shifted Gaussians share the same poles, which reduces the terms in the next step. Then, they determine the coefficients $a_l$ by minimizing the following $l^{\infty}$ error
\begin{equation}
\label{eq:minlinfapprox}
\max_{x_k} \left|\frac{1}{\sqrt{4\pi}}e^{-x_k^2/4} - R(x_k)\right|
\end{equation}
using a set of points $x_k \in [-30,30]$. In Table~1 of \cite{HBMW2015}  the values of the coefficients $a_l$ and $\mu$ are given. The corresponding approximation error is less than $7 \times 10^{-13}$. 

We now propose a strategy to reduce the approximation error by using a different approach to determine the coefficients. First, we observe that the function $\psi_1(x)$ is a symmetric function, i.e.~$\psi_1(x)=\psi_1(-x)$. Therefore, the same property should hold true for the corresponding approximation $R(x)$, which is not satisfied in \cite{HBMW2015}. We remark that $R(x)$ is symmetric if $a_l = \overline{a_{-l}}$ for all $l$. To enforce that this property holds also numerically, we rewrite \eqref{eq:ratapproxgaush1} as
\begin{equation}
\label{eq:A(x,mu)}
R(x) = \frac{a_0\mu }{x^2 + \mu^2} + \sum_{l = 1}^L \frac{2\mu \re(a_l)
(\mu^2 + l^2 + x^2) + 2l\im(a_l)(\mu^2 + l^2 - x^2)}{x^4 + 2(\mu^2 - l^2)x^2 + (\mu^2 + l^2)^2}.
\end{equation}

Then, to find the coefficients $a_l$ we minimize the following $l^2$ error
\begin{equation}
\label{eq:minl2approx}
\sum_{k=1}^K \left(\frac{1}{\sqrt{4\pi}}e^{-x_k^2/4} - R(x_k)\right)^2
\end{equation}
on $K$ points $x_1, \dots , x_K$. 
Since $R(x)$ is a linear combination of the coefficients, the approximation can be written as
\[
R(x_k) = G(x_k,\mu,L)y \approx \psi_1(x_k),
\]
where $G(x_k,\mu,L) \in \mathbb{R}^{1 \times (2L+1)}$ and $y \in \mathbb{R}^{2L+1}$ is of the form
\[
y = [a_0, \re(a_1), \re(a_2), \dots, \re(a_L), \im(a_1), \dots, \im(a_L)]^T.
\]
The reason why we chose to minimize \eqref{eq:minl2approx} instead of \eqref{eq:minlinfapprox} is because it can be solved easier and faster. To minimize \eqref{eq:minl2approx} we compute the least square solution of the corresponding linear system, where the points $x_k$ are calculated iteratively. We start with $x_1 = 0$ and for selecting the next point $x_{k+1}$ we use the same strategy that is used for minimizing the error in interpolation with Leja points~\cite{Leja}. For $L = 24$ we obtained the coefficients listed in Table~\ref{tab:coef_al}. The coefficient $\mu$ is determined such that a high accuracy is obtained. For this choice the error in the maximum norm is less than $8\times 10^{-15}$.

\subsubsection*{Step 3}

The third step is the combination of steps 1 and 2:
\begin{equation}
\label{eq:RexiNotSimplified}
e^{ix} \approx \sum_{m = -M}^M b_m \sum_{l = -L}^L \re\left(
\frac{ha_l}{ix + h(\mu + i(m+l))}\right).
\end{equation}

For computational efficiency equation \eqref{eq:RexiNotSimplified} should be rewritten as a single sum. This can be done in different ways. A possibility is to split up the coefficients $b_m$ into their real and imaginary parts and pull them inside the approximation of the Gaussian function. This is done in \cite{SPHW2017}. Let  $n = m + l, N = M + L, \alpha_n = h(\mu + in)$ and set
\[
\beta_n^{Re} = h \sum_{k = L_1(n)}^{L_2(n)} a_k\re(b_{n-k})
\]
and
\[
\beta_n^{Im} = h \sum_{k = L_1(n)}^{L_2(n)} a_k\im(b_{n-k}),
\] 
where $L_1(n) = \max(-L,n-M)$, $L_2(n) = \min(L,n+M)$. This leads to the following form 

\begin{equation}
\label{eq:originalRexi}
\rexi(ix) := \sum_{n = -N}^{N} \re\left(\frac{\beta_n^{Re}}{ix + \alpha_n}\right) + i\re\left(\frac{\beta_n^{Im}}{ix + \alpha_n}\right),
\end{equation} 
which serves as a definition of the numerical approximation $\rexi(ix)$ to $e^{ix}$.

Another possibility to rewrite \eqref{eq:RexiNotSimplified} is to compute the real part of the approximation of the Gaussian function. In contrast to the previous reformulation, we make now explicit use of the fact that $x$ has to be real. We set
\[
c_{1,n} = h \sum_{k = L_1(n)}^{L_2(n)} \re(a_k)b_{n-k}
\] 
and
\[
c_{2,n} = h \sum_{k = L_1(n)}^{L_2(n)} \im(a_k)b_{n-k},
\] 
and obtain an equivalent form of \eqref{eq:originalRexi} which we call $\rexii$, thus
\begin{equation}
\label{eq:modifiedRexi}
\rexii(ix) := \sum_{n = -N}^{N} \frac{c_{1,n}h\mu + c_{2,n}(x+hn)}{(\alpha_{-n}-ix)(\alpha_n +ix)}.
\end{equation}
Recall that $N = M+L$ and $\alpha_n = h(\mu+in)$.  In the scalar case both simplifications take roughly the same computational effort and the results are equivalent for real $x$. This, however, is not true in the matrix case (as we will see in the next section).

We now explain how \eqref{eq:modifiedRexi} can be used to compute $e^{\tau A}$ for a given square matrix with purely imaginary eigenvalues. Assume that $A$ is diagonalizable by a matrix $V$, i.e.,
\begin{equation}
\label{eq:AisVEV}
A = VEV^{-1}
\end{equation}
where $E$ is a diagonal matrix. Then it follows that $e^{\tau A} = Ve^{\tau E}V^{-1}$. The diagonal entries of $e^{\tau E}$ can be computed componentwise by REXII
\begin{equation}
\label{eq:exponential_of_E}
e^{\tau E} = 
  \begin{pmatrix}
    \ddots &  &   \\
     & e^{i\lambda_j\tau}  &\\
     & & \ddots
  \end{pmatrix}
\approx
  \begin{pmatrix}
    \ddots &  &   \\
     & \rexii(i\lambda_j\tau)  &\\
     & & \ddots
  \end{pmatrix}
=: \rexii(\tau E),
\end{equation}
where $i\lambda_j$ are the eigenvalues of $A$.  Finally we obtain
\begin{equation}
\label{eq:REXI_VEV_DECOMP}
e^{\tau A} \approx V\rexii(\tau E)V^{-1}.
\end{equation}
This is a well known technique to extend scalar functions to matrices, where the scalar functions are applied to the spectrum of the matrix. Therefore, $\tau A$ can be substituted for $ix$ in \eqref{eq:modifiedRexi} since $A$ is a matrix with purely imaginary eigenvalues. 

\subsection{The new scheme $\rexii$ for matrices}
\label{sec:REXIMATRIXNEW}

Let $A$ be a matrix with purely imaginary eigenvalues. As explained before, we can substitute $\tau A$ for $ix$ in \eqref{eq:modifiedRexi}. This yields 
\[
e^{\tau A} \approx \sum_{n = -N}^N (c_{1,n}h\mu I + c_{2,n}(-i\tau A+hnI))(\alpha_{-n}I-\tau A)^{-1}(\alpha_n I+\tau A)^{-1}.
\]

This scheme can be made more efficient by defining $C_{1,n} = (c_{1,n}h\mu + c_{2,n}hn)$ and $C_{2,n} = ic_{2,n}$. As
\begin{align*}
(C_{1,n}I - C_{2,n}\tau A)&(\alpha_{-n}I - \tau A)^{-1} = \\
 &C_{2,n}\left(\left(\frac{C_{1,n}}{C_{2,n}} - \alpha_{-n}\right)(\alpha_{-n}I - \tau A)^{-1} +I \right)
\end{align*} 
we obtain the following extension of $\rexii$ for matrices with purely imaginary eigenvalues
\begin{equation}
\label{eq:REXI_Modified_matrix}
\begin{aligned}
\rexii(\tau A) := \sum_{n = -N}^N C_{2,n}\left(\left(\frac{C_{1,n}}{C_{2,n}} - \alpha_{-n}\right)(\alpha_{-n}I - \tau A)^{-1} +I \right)(\alpha_n I+\tau A)^{-1}. \\
\end{aligned}
\end{equation}
Recall that $N = M+L$ where $L=24$ and $\alpha_n = h(\mu + in)$.  We are now in a position to show how for $\rexii$ the accuracy of the scalar case translates to the matrix case. This is the content of Theorem~1. Its proof is a consequence of \eqref{eq:AisVEV} and \eqref{eq:exponential_of_E}.

\begin{Teo}
\label{teo:rexi}
Let $A$ be a square matrix and suppose that $\sigma (A) \subset i\mathbb{R}$. If $A$ is diagonalizable, i.e., $A = VEV^{-1}$ with $E = \diag(i\lambda_j)$, then
\begin{equation}
\label{eq:teorem}
\|\rexii(\tau A) - e^{\tau A}\|_{\infty} \leq  \cond(V)\cdot \max_{\lambda_j \in \sigma(A)}|\rexii(i\tau\lambda_j) - e^{i\tau \lambda_j}|
\end{equation}
where $\cond(V) = \|V\|_{\infty} \cdot \|V^{-1}\|_{\infty}$. 
\end{Teo}

Regarding the error in the matrix case, by Theorem~\ref{teo:rexi} we have to estimate $\rexii$ in the scalar case for every eigenvalue $i\tau \lambda_j$ of $\tau A$. In the scalar case $\rexii$ is accurate if condition \eqref{eq:Mformula} holds for every $i\tau \lambda_j$, namely $|\tau \lambda_j| \leq (M-11)h$. This implies in the matrix case, that if we choose $M$ and $h$ such that
\begin{equation}
\label{eq:matrixAccuracyBound}
\tau \rho(A) \leq (M-11)h
\end{equation}
holds true, where $\rho(A)$ is the spectral radius of $A$, we are guaranteed to obtain results close to machine precision. In practice, since $\rexii$ in general is applied to a vector, also the error from solving the linear systems has to be taken into account.

\begin{rmk}
From \eqref{eq:matrixAccuracyBound} we observe that it might be possible to reduce the amount of work by using a shifted matrix, i.e.~$A' = A-\nu I$.  If $\rho(A') < \rho(A)$ then less terms are needed to obtain an accurate approximation of $e^{A'}$. Note that the matrix $e^A$ can easily be recovered from $e^{A'}$ as follows
\[
e^A = e^{A -\nu I + \nu I} = e^\nu e^{A - \nu I} = e^\nu e^{A'}.
\]
Now, since in our case all eigenvalues are purely imaginary, we have $i\lambda_j \in i[\zeta_1,\zeta_2]$. Thus to reduce $\rho(A)$ we use the shift $\nu = i\frac{\zeta_1 + \zeta_2}{2}$. This implies that $\rho(A') = \zeta_1-\nu$ and $i\lambda'_j \in i[\zeta_1-\nu,\zeta_2-\nu]$. For skew symmetric real matrices, it is not necessary to perform a shift since all eigenvalues arise in complex conjugate pairs, and therefore $\zeta_1 = -\zeta_2$.
\end{rmk}

\begin{rmk}
If the matrix $A = VEV^{-1}$ is skew Hermitian, then the matrix $V$ is unitary. If moreover, the error in \eqref{eq:teorem} is estimated in 2-norm, then $\cond(V) = 1$ and therefore the error in the matrix case is the same as in the scalar case for the eigenvalues of $A$.
\end{rmk}

\begin{rmk}
Since $\rexii$ in general is used to evaluate the action of the matrix exponential applied to a vector $f_0$, the main cost is to solve two linear systems for each summation term. We are able to reduce the cost if the entries of $A$ and $f_0$ are real. In particular, we observe that $\overline{c_{1,n}} = c_{1,-n}$, $\overline{c_{2,n}} = -c_{2,-n}$ and $\alpha_n = \overline{\alpha_{-n}}$. Therefore, we obtain
\begin{equation}
\label{eq:modifiedRexiMatrixReducedSum}
e^{\tau A} \approx \re\left( \sum_{n = 0}^N \Gamma_nC_{2,n}\left(\left(\frac{C_{1,n}}{C_{2,n}} - \alpha_{-n}\right)(\alpha_{-n}I - \tau A)^{-1} +I \right)(\alpha_n I+\tau A)^{-1}\right),
\end{equation}
where $\Gamma_0 = 1$ and $\Gamma_n = 2 \enspace \text{for} \enspace 1 \leq n \leq N$.
\end{rmk}

\subsection{The original $\rexi$ scheme and the differences to $\rexii$}
\label{sec:REXIMATRIXOLD}
The original REXI scheme was developed for matrices $A$ with real entries. It is based on~\eqref{eq:originalRexi} where $\tau A$ is substituted for $ix$. A further simplification comes from the fact that $e^{\tau A}$ is real which suggests to neglect the imaginary part of~\eqref{eq:originalRexi}. The scheme is thus defined as follows
\begin{equation}
\label{eq:originalREXImatrix}
\rexi(\tau A) = \sum_{n = -N}^N \re\bigl(\beta_n^{Re}(\tau A + \alpha_n I)^{-1}\bigr) \bigr).
\end{equation}

At first look this approximation to $e^{\tau A}$ might be the preferred one since if $\rexi$ is applied to a real vector $f_0$, only one linear system has to be solved for each summation term. In contrast, $\rexii$ has two linear systems to solve for each term. But $\rexi$ has some drawbacks with respect to $\rexii$.
\begin{itemize}
\item First, if $A$ is real with purely imaginary eigenvalues, then $V$ has to be complex and therefore in \eqref{eq:REXI_VEV_DECOMP} the matrix $V$ can not be pulled inside the real part of the approximation. Thus, Theorem~\ref{teo:rexi} does not hold for $\rexi$ and condition \eqref{eq:matrixAccuracyBound} does not apply. To improve the accuracy, $M$ has to be increased. The rate of convergence of $\rexi$ can be slow and an extremely large value for $M$ might be  needed if a stringent tolerance is prescribed, as we will show in the numerical experiments.  Therefore, the main advantage of $\rexii$ compared to $\rexi$ is relation \eqref{eq:matrixAccuracyBound} that allows us to choose $h$ and $M$ in an appropriate manner to produce the same high accuracy as in the scalar case. The cost of this is that two linear systems have to be solved, and thus the sequential part of the scheme doubles. Despite the increased computational effort, however, our scheme is still significantly faster since we can choose a much smaller $M$.

\item Second, reducing the sum from $2N+1$ terms to $N+1$ terms, as is done in \eqref{eq:modifiedRexiMatrixReducedSum}, is not feasible for the original REXI scheme since the coefficients $a_l$ in \cite{HBMW2015} are not exactly equal to $\overline{a_{-l}}$. Doing this, as in \cite{SPHW2017}, result in a reduction of accuracy in the approximation of the Gaussian functions (from $7 \times 10^{-13}$ to $4\times 10^{-8}$) and since $\rexi$ is at most as accurate as the approximation of the Gaussian function, the overall accuracy is significantly reduced.
\end{itemize}

\begin{rmk}
If $A = iB$, where $B$ is a real diagonalizable matrix, then the transformation matrix $V$ is real. If moreover $f_0$ is real, the use of \eqref{eq:originalRexi} for approximating $e^{\tau A}$ is justified, and we end up with the following scheme:
\begin{equation}
\label{eq:REXIE}
e^{\tau A}f_0 \approx \sum_{n=-N}^{n=N} \re\left(\beta_n^{\re}(\tau A + \alpha_nI)^{-1}f_0\right) + i\re\left(\beta_n^{\im}(\tau A + \alpha_nI)^{-1}f_0\right)
\end{equation}
We call this scheme $\rexi$ Extended (REXIE). Note that this formulation is more efficient than $\rexii$ as only one linear system has to be solved. Matrices of the above form are, for example, purely imaginary skew Hermitian matrices.
\end{rmk}

\section{Numerical examples}
\label{sec:NumericalExamples}

In this section we provide numerical examples that confirm the theoretical considerations laid out in the previous section. We start with the scalar case and then advance to the matrix case. We always use the conjugate symmetric coefficients given in Table~\ref{tab:coef_al} for the implementation of $\rexi$. This is a slight modification of the algorithm in \cite{SPHW2017} but reduces cost and improves accuracy, as explained in the previous section.

\subsection{The scalar case}
\label{sec:numericalExamplesScalar}

Note that in the scalar case, the absolut error is the same as the relative error, since $|e^{ix}| = 1$. All the results in this sections are calculated sequentially using GNU~Octave.

\begin{figure}[H]
  \centering
  \subfloat[][]{\includegraphics[width=0.47\textwidth]{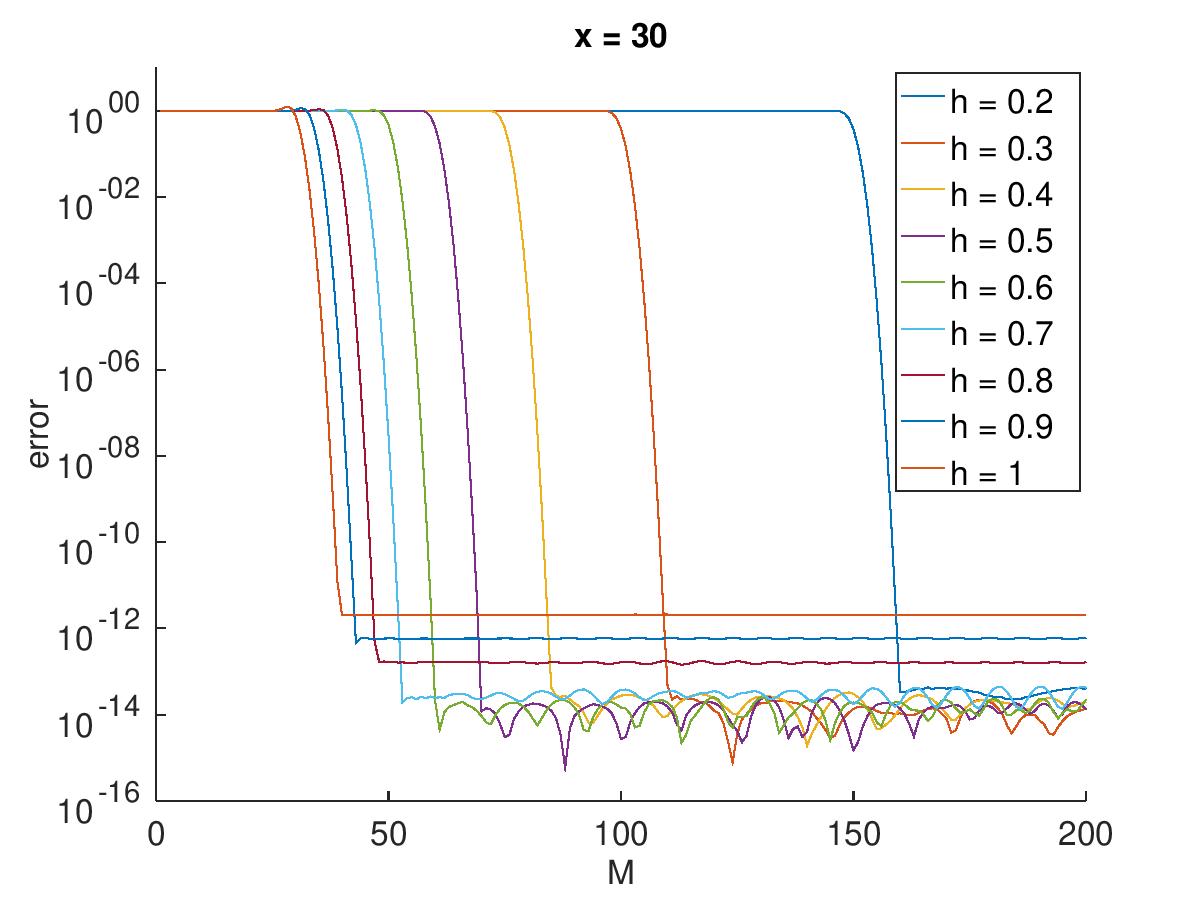}}%
  \qquad
  \subfloat[][]{\includegraphics[width=0.47\textwidth]{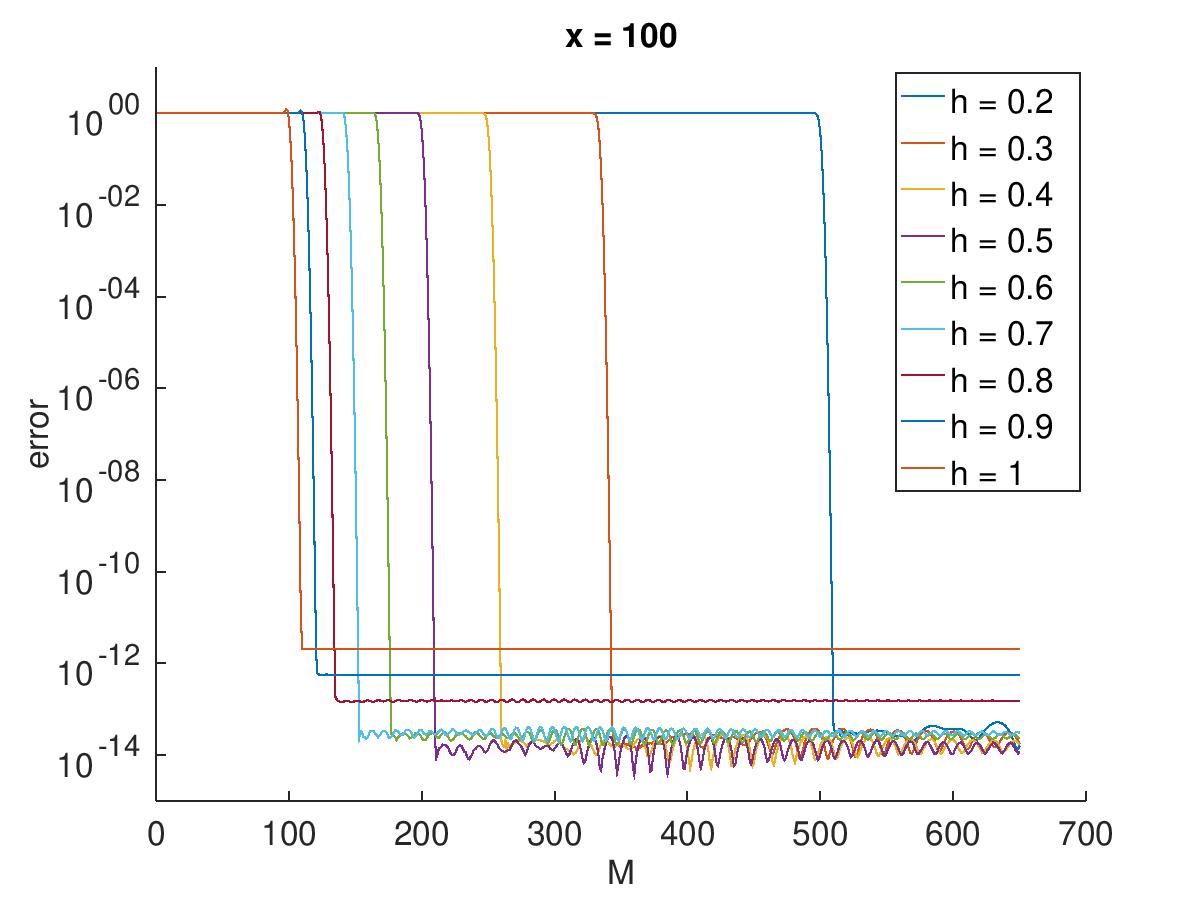}}%
  \caption{\label{fig:exLB}Error of the approximation \eqref{eq:originalRexi} to $e^{30i}$ and $e^{100i}$ for different choices of $h$ as a function of the parameter $M$.}%
\end{figure}

In order to study the approximation, we fix $x$ and plot the error as a function of $h$ and $M$. The results are shown in Figure~\ref{fig:exLB}. We observe that after a certain value of $M$ the error drops immediately to a value close to machine precision. The point at which this happens is well predicted by the bound $|x| \leq (M-11)h$. Furthermore, we observe that we obtain the best results for $h \in [0.3,0.6]$.  Thus it is not recommended to choose $h$ too large or too small. A small value of $h$ leads to a big value of $M$,  which increases the computational cost. A large value of $h$ results in reduced accuracy.

We have also investigated the error as a function of $x$, where we fixed $h$ and let $M$ to be the minimum admissible values given from \eqref{eq:Mformula}. We observe that for this numerical test the accuracy is always close to machine precision.

\subsection{The matrix case}
\label{sec:numericalExamplesMatrix}
\begin{figure}[H]
  \centering
  \subfloat[][]{\includegraphics[width=0.47\textwidth]{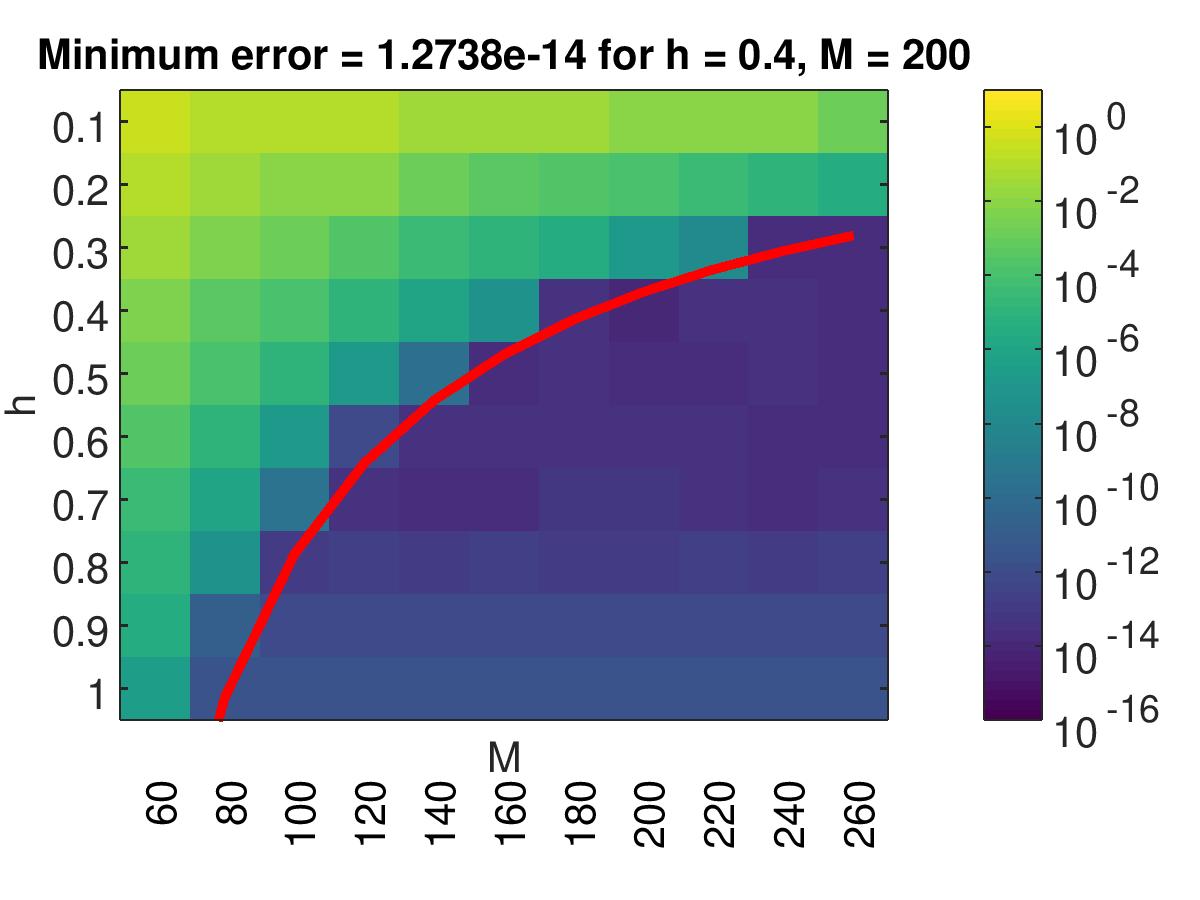}}%
  \qquad
  \subfloat[][]{\includegraphics[width=0.47\textwidth]{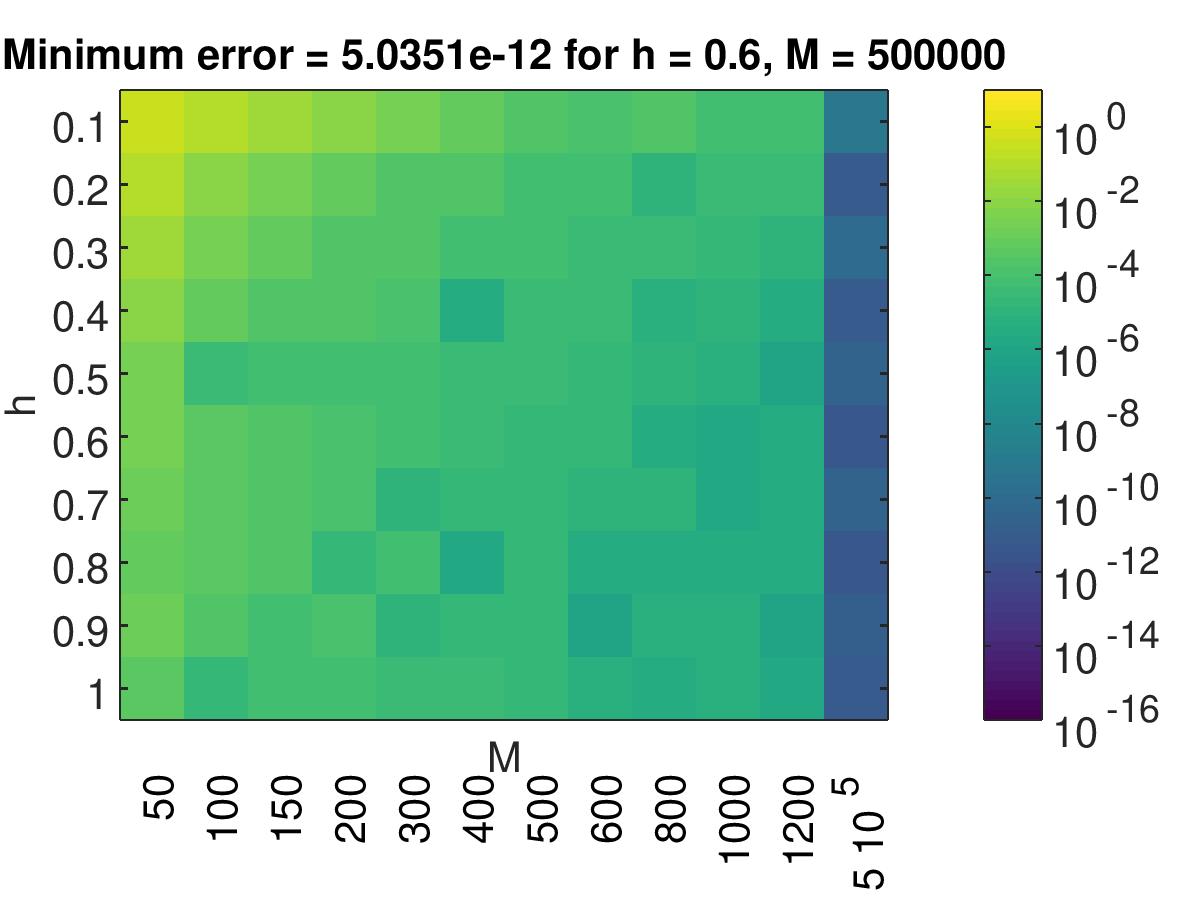}}%
  \qquad
  \subfloat[][]{\includegraphics[width=0.47\textwidth]{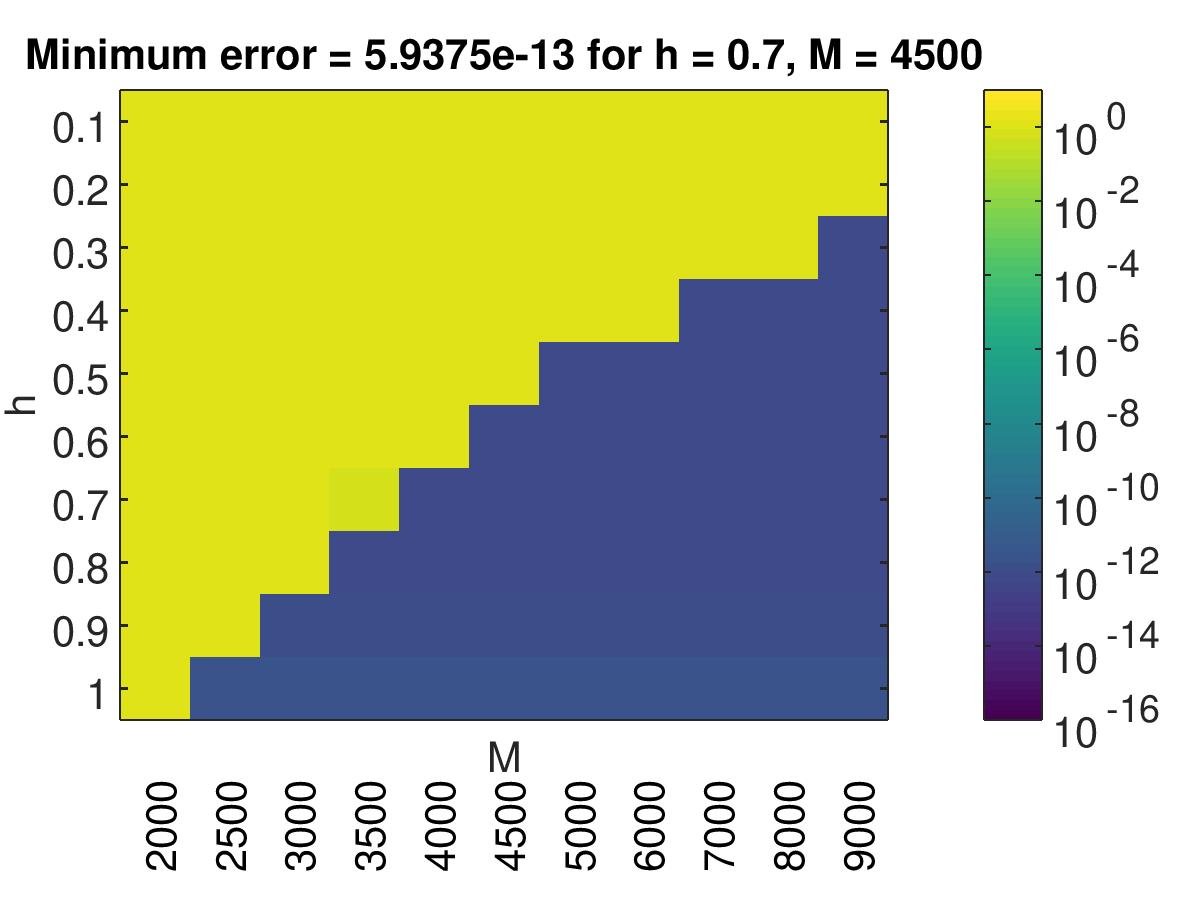}}%
  \caption{\label{fig:matrixRexiExample} Relative error of the approximation of $e^{A_1}f_0$ with REXII in (a) and with the original REXI in (b). The red curve represents the bound \eqref{eq:matrixAccuracyBound}, namely $\rho(A) = (M-11)h$. For (b) we can see that this bound does not provide a good prediction. Furthermore, we need a very large $M$ to produce accurate results compared to the scheme in (a). The error of the approximation of $e^{A_2}f_0$ with REXIE \eqref{eq:REXIE} is given in (c). }%
\end{figure}

In this section, we analyze the behavior of the different algorithms depending on the parameters $h$ and $M$ for two test matrices $A_1$ and $A_2$.
\begin{itemize}
\itemsep0em 
\item The matrix $A_1$ is the second order finite difference approximation of the advection operator $\partial_x$ with periodic boundary conditions in the spatial domain $[a,b] = [0,1]$. The discretization step is $\frac{b-a}{n-1} = 1/70$. Therefore, $A_1$ is a skew symmetric matrix with eigenvalues $\lambda_j \in i[-70,70]$. This matrix has complex eigenvectors, thus we apply $\rexii$ \eqref{eq:modifiedRexiMatrixReducedSum}. With this matrix $\rexii$ has to solve two linear systems.

\item The matrix $A_2$ is the second order finite difference approximation of the free Schrö\-din\-ger operator $i\partial_{xx}$ in the spatial domain $[-1,1]$ with periodic boundary conditions. The discretization step is $1/35$. Thus $A_2$ is a purely imaginary skew Hermitian matrix with eigenvalues $\lambda_j \in i[-4900,0]$. Therefore, it is convenient to apply a shift of $\nu = -2450i$. This matrix has real eigenvectors, thus we apply  REXIE \eqref{eq:REXIE}. With this matrix REXIE has to solve only one linear system per summation term.
 \end{itemize}
As vector $f_0$ in both cases we used the discretization of $(2+\cos(2\pi x))^{-1}$. The results in this section are computed in GNU Octave. We measure the relative error in the $l^2$ norm

\[
    \text{err} = \frac{\|\rexii(A)f_0 - \text{{\fontfamily{qcr}\selectfont expm}}(A)f_0\|_2}{\|\text{{\fontfamily{qcr}\selectfont expm}}(A)f_0\|_2},
\]
where {\fontfamily{qcr}\selectfont expm} is a Pad\'e approximation of the exponential matrix. The results are shown in Figure~\ref{fig:matrixRexiExample}. We clearly see that $\rexii$ is much more accurate for the same $M$ compared to the original $\rexi$ scheme. In addition, the bound \eqref{eq:matrixAccuracyBound} predicts the behavior of $\rexii$ very well, which is not the case for the $\rexi$ scheme.

\section{Linear rotating shallow water equations}
\label{sec:numericalExampleLSWE}
In this section we apply our proposed REXII scheme and the original REXI scheme to the linear rotating shallow water equations (LRSW) \cite{LSWE}. This is the same problem that has been investigated in \cite{SPHW2017}. The LRSW are stated as follows
\begin{align*}
\partial_t f = Af, \quad f(0) = f_0,
\end{align*}
where the linear operator $A$ is defined as follows
\[
A=
\begin{pmatrix}
0 			& -\partial_x   & -\partial_y \\
-\partial_x & 0				& 1 \\
-\partial_y & -1 			& 0\\ 
\end{pmatrix}.
\]
The sought after function is $f = (\eta, u ,v)^T$, where $\eta$ is the displacement of the surface height, $u$ the velocity in the $x$ direction and $v$ the velocity in the $y$ direction. The simulation domain is the bi-periodic unit square $[0,1]^2$, thus periodic in the $x$ direction, $f(\tau,0,y) = f(\tau,1,y)$, and periodic in the $y$ direction, $f(\tau,x,0) = f(\tau,x,1)$. The grid resolution is $D \times D = 128 \times 128$. To apply REXII \eqref{eq:modifiedRexiMatrixReducedSum}, we have to solve for each term two linear systems and sum up these two calculated solutions. More specifically we have to compute for $0 \leq n \leq N$

\begin{align*}
(A + \alpha_n I)g_{1,n} &= f_0 \\
(\alpha_{-n} I - A)g_{2,n} &= g_{1,n} \\
g_{3,n} &= \Gamma_nC_{2,n}(g_{1,n} + (C_{1,n}/C_{2,n} - \alpha_{-n})g_{2,n}) \\
\intertext{and}
e^{A}f_0 &= \sum_{n = 0}^N \re( g_{3,n}).
\end{align*}
To solve the linear systems, the following strategy is applied (see also \cite{SPHW2017}). Taking the second and the third component of $(A + \alpha_n I)g_{1,n} = f_0$ leads to the following equation for the velocities
\[
\begin{pmatrix}
\alpha_n	& 1   		 \\
-1 			& \alpha_n	 \\
\end{pmatrix}
\begin{pmatrix}
u_{1,n}	\\
v_{1,n} 	\\
\end{pmatrix}
= 
\begin{pmatrix}
u_0			\\
v_0 			\\
\end{pmatrix}
+
\nabla \eta_{1,n}.
\]
Inverting the $2\times2$ linear system yields
\begin{equation}
\label{eq:lswVelocities}
\begin{pmatrix}
u_{1,n}	\\
v_{1,n} 	\\
\end{pmatrix}
= 
\frac{1}{\kappa_n}
\begin{pmatrix}
\alpha_n	& -1   		 \\
1 			& \alpha_n	 \\
\end{pmatrix}
\left(
\begin{pmatrix}
u_0		   \\
v_0 			\\
\end{pmatrix}
+
\nabla \eta_{1,n}
\right),
\end{equation}
where $\kappa_n = 1 + \alpha_n^2$. The velocities can then be calculated directly if $\eta_{1,n}$ is available. From the first component of $(A + \alpha_n I)g_{1,n} = f_0$ we obtain
\begin{equation}
\label{eq:lswEtaIndirect}
\alpha_n \eta_{1,n}- \nabla \cdot 
\begin{pmatrix}
u_{1,n} \\
v_{1,n} \\
\end{pmatrix}
= \eta_0.
\end{equation}
Thus, if we plug \eqref{eq:lswVelocities} into \eqref{eq:lswEtaIndirect} we obtain
\begin{equation}
\label{eq:lswEta}
\Delta \eta_{1,n} - \kappa_n \eta_{1,n} = r_{0,n},
\end{equation}
where 
\begin{align*}
r_{0,n} &= -\frac{\kappa_n}{\alpha_n} \eta_0 + \frac{1}{\alpha_n}\zeta_0 - \delta_0 \\
\delta &= u_x + v_y \\
\zeta &= v_x - u_y.
\end{align*}
Therefore, we first solve the Helmholtz problem \eqref{eq:lswEta} for $\eta$. Then we obtain the velocities by plugging the gradient of $\eta$ into \eqref{eq:lswVelocities}. To efficiently solve the linear system, all computations will be conducted in Fourier space.

The parameter $M$ of REXI depends on the spectral radius of $A$. We obtain
\begin{equation}
\label{eq:lswRoh}
    \rho(A) = \sqrt{2 \pi^2 D^2+1} \approx \sqrt{2} \pi D.
\end{equation}

We will use different initial conditions to illustrate this example. We compare REXII with REXI and the explicit Runge--Kutta time stepping method of  order 4 (RK4). For the results in the following test scenarios, we use the maximum norm:
\[
    \text{err} = \max_{r,s}|f_{\text{num}}(\tau,x_r,y_s) - f_{\text{ref}}(\tau,x_r,y_s)|,
\]
where $f_{\text{num}}(\tau)$ is the numerical approximation and $f_{\text{ref}}(\tau)$ is a reference solution. For REXI and REXII we perform only one single time step of size~$\tau$, and $(x_r,y_s)$ are the grid points of the domain. Since it is not clear for REXI how to choose $M$ and $h$, we fix $h = 0.2$ as is done in \cite{SPHW2017} and vary $M$. Before presenting the numerical results we discuss the parallel implementations for both CPU and GPU based systems.

\subsection{Implementation}
\label{sec:Implementation}

The calculations in this section are conducted on a single GPU (NVIDIA V100) and separately on a CPU (a dual socket Intel Xeon Gold 5118 server with a total of 32 cores). We also performed the calculations on a NVIDIA TitanV. Since the performance of both NVIDIA cards is quite similar we only report the V100 results here. The code is written in C++ and CUDA 10.0 is used to program the GPU. For the GPU code we use CUFFT \cite{CUFFT} to perform FFTs and to do the reduction sum at the end we use the CUB library \cite{CUB}. For the CPU code we use the FFTW library~\cite{FFTW}. The code for the CPU is implemented sequentially over the summation index $n$, but each iteration is parallelized with OpenMP~\cite{OpenMP}. Initially, we used the BLAS implementation found in Intel MKL~\cite{MKL} for the CPU code. However, due to the possibility to aggregate many operations, which is a big advantage for memory bound problems, our OpenMP implementation is actually significantly faster. The GPU code parallelizes over the sum in addition to the spatial parallelization. This is done to exploit the massively parallel architecture of modern GPUs.

\subsubsection*{CPU Implementation}

The implementation for the CPU is sequential over the summation range $0:N$ and is therefore performed exactly how it is described in the previous section. We denote by $\hat{\eta}_{1,n}$ the solution of the variable $\eta$ of the first system and $\hat{\eta}_{2,n}$ the solution of the variable $\eta$ of the second system where $n$ is the index in the sum of REXI. Similar notation holds for the other variables. To perform the calculation of these variables, we use OpenMP with 32 threads.

\begin{algorithm}[H]
  \caption{CPU implementation}
  \begin{algorithmic}
  	\STATE $[\hat{\eta}_0, \hat{u}_0, \hat{v}_0 ]= \text{FFTW}(\eta_0,u_0,v_0)$
  	\STATE $\text{calculate } \hat{\delta}_0, \hat{\zeta}_0 $
  	\STATE $\text{calculate all }C_{1,n}, C_{2,n}$
    \FOR{$n = -N:0$} 
    	\STATE $\text{solve the first system, calculate }\hat{\eta}_{1,n}, \hat{u}_{1,n}, \hat{v}_{1,n}$
		\STATE $\text{prepare to solve the second system, calculate } \hat{\delta}_{1,n}, \hat{\zeta}_{1,n}$
		\STATE $\text{solve the second system, calculate } \hat{\eta}_{2,n}, \hat{u}_{2,n}, \hat{v}_{2,n}$
		\STATE $\text{reduce: } \hat{\eta}= \hat{\eta} + (C_{2,n}\hat{\eta}_{1,n} + C_{1,n}\hat{\eta}_{2,n}), \text{ similarly for } \hat{u}, \hat{v}$
    \ENDFOR
    \STATE $[\eta, u, v] = \re(\text{IFFTW}(\hat{\eta}, \hat{u}, \hat{v}))$
    \label{alg:CPUimplementation}
  \end{algorithmic}
\end{algorithm}
The achieved memory bandwidth given in Table~\ref{tab:memoryBandwidth} confirms that the implementation performs as expected.

\subsubsection*{GPU Implementation}

For the GPU implementation we compute the linear systems in parallel. To facilitate the computation of $\hat{\eta}$ we create a matrix of size $(N+1) \times D^2$, where each row corresponds to $\hat{\eta}_n$ for $n = 0 \ldots N$. We denote this matrix by E1. In a similar way we create two other matrices for $\hat{u}$ and $\hat{v}$, which we call U1 and V1, respectively. We store these matrices in column major format, such that we can apply the fast reduction library CUB to perform the sum over $n$. Therefore, for REXI we need three such matrices. For REXII we need six of these matrices because we have to solve two different linear systems. We denote these matrices by E1, E2, U1, U2, V1, V2. These matrices are very memory intensive. Each entry requires 16 Bytes of memory, since we are using complex double precision arithmetic. Therefore, for REXII for example, we  need approximately $6 \cdot 16 \cdot D^3 \cdot \sqrt{2} \pi \tau /h$ bytes of memory to store these six matrices, which corresponds to approximately 90~GBytes for $\tau = 50$, $h = 0.5$ and $D = 128$. The GPU we are working with has only 12~GBytes of memory. Therefore, we have to divide the total amount of required memory to store the variables in $S$ parts. Thus we are calculating at each step only a part of E1, namely $\text{E1}_s = \text{E1}(s\bar{N}:(s+1)\bar{N}, :)$ where $s = 0\ldots S-1$ and $\bar{N} =  (N+1)/S $.

There are six kernel functions required. They carry out the following tasks:
\begin{itemize}
\itemsep0em 
\item computation of E1$_s$
\item computation of U1$_s$ and V1$_s$
\item computation of $\delta$ and $\zeta$ to solve the second linear system
\item computation of E2$_s$
\item computation of the final data
\item reduction via CUB library
\end{itemize}

The corresponding algorithm is summarized in Algorithm \ref{alg:GPUimplementation} and the overall achieved performance is listed in Table ~\ref{tab:memoryBandwidth}. We clearly see that the GPU implementation outperforms the CPU implementation by a significant margin.

\begin{algorithm}[H]
  \caption{\label{alg:GPUimplementation} GPU implementation}
  \begin{algorithmic}
  	\STATE $\text{FFT of initial data}$
  	\STATE $\text{divide the required memory in }S\text{ parts}$
    \FOR{$s = 0:S-1$} 
    	\STATE $\text{solve the first linear system:}$
    	\STATE $\text{E1}_s = \text{computeE1}(\hat{\eta}_0, \hat{\delta}_0, \hat{\zeta}_0, ...);$
		\STATE $[\text{U1}_s, \text{V1}_s] =  \text{computeU1V1}(\hat{u}_0, \hat{v}_0, \text{E1}_s,  ...);$
		\STATE $\text{solve the second linear system and reduce:}$
		\STATE $[\text{U2}_s, \text{V2}_s] = \text{computeDelta1Zeta1}(\text{U1}_s, \text{V1}_s);$
		\STATE $\text{E2}_s = \text{computeE2}(\text{E1}_s, \text{U2}_s, \text{V2}_s, ...);$
		\STATE $[\text{E2}_s, \text{U2}_s, \text{V2}_s] = \text{finalize}(\text{E1}_s, \text{U1}_s, \text{V1}_s, \text{E2}_s, ...	);$
		\STATE $\text{reduction via CUB}$
    \ENDFOR
    \STATE $\text{inverse FFT to obtain solution}$
  \end{algorithmic}
\end{algorithm}

\begin{table}
\begin{center}
\begin{tabular}{l|c|c|c}
 & GB/s & write and read & expected factor \\ \hline
REXI GPU & 480 & 7 & \multirow{2}{*}{17.1} \\ 
REXI CPU & 80 & 20 &  \\ 
\hline 
REXII GPU & 370 & 22 & \multirow{2}{*}{9.7} \\ 
REXII CPU & 80 & 46 &  \\ 
\hline 
RK4 GPU & 200 & 42 & \multirow{2}{*}{4.7} \\ 
RK4 CPU & 60 & 59 &  \\ 
\end{tabular} 
\caption{ \label{tab:memoryBandwidth} The memory bandwidth and the number of memory write and read operations for each iteration of the algorithm is listed. Form these numbers we can deduce a expected factor of speedup on the GPU compared to the CPU. On the GPU, we need less memory operations because the implementation computes large blocks of data at once.}
\end{center}
\end{table}

\subsection{Wave scenario 1}
\label{sec:ws1}
The following initial conditions are used for wave scenario 1:
\begin{equation}
\label{eq:WAVESCENARIO1}
\begin{aligned}
\eta(0,x,y) &= \sin(4\pi x)\cos(2\pi y) -\tfrac{1}{5}\cos(4\pi x)\sin(4 \pi y) \\
u(0,x,y) &= \cos(8\pi x)\cos(2\pi y) \\
v(0,x,y) &= \cos(4\pi x)\cos(4\pi y)
\end{aligned}
\end{equation}
This is the same problem as considered in \cite{SPHW2017}. Since we solve the problem in Fourier space, these initial functions are extremely convenient. They are exactly representable in Fourier space with very few terms. This gives us a big advantage when we have to choose the parameter $M$ of REXII. In \eqref{eq:lswRoh} we can consider $D = 6$, since all higher modes do not contribute to the solution.
From \eqref{eq:matrixAccuracyBound} we deduce that $M$ is
\[
M = \left\lceil \frac{\sqrt{2}\cdot \pi\cdot 6\cdot \tau}{h} \right\rceil + 11,
\]
where $\tau$ is the final time, since we are performing only one time step. We thus include this problem primarly to provide a comparison to the results obtained in \cite{SPHW2017}. We will conduct an investigation with more realistic initial values in the subsequent sections.

\begin{table}
\begin{center}
\begin{tabular}{ c| r@{\hskip 4pt}r@{\hskip 1pt}r | r@{\hskip 1pt}l | r@{\hskip 1pt}l |r@{\hskip 1pt} l@{\hskip 2pt}l |r@{\hskip 1pt} l@{\hskip 2pt}l }
\multicolumn{12}{ c }{Final time $\tau$ = 1} \\
\hline
\multirow{2}{*}{Method} & \multicolumn{3}{c}{\multirow{2}{*}{$h,M$ \hfill / \hfill time steps TS}} & \multicolumn{4}{c}{Error} & \multicolumn{6}{c}{Time}\\
& & & & \multicolumn{2}{c|}{CPU} & \multicolumn{2}{c|}{GPU} & \multicolumn{3}{c|}{CPU} & \multicolumn{3}{c}{GPU} \\ \hline
\multirow{3}{*}{REXI} 
 & $h$ = 0.2, $M$ =& & 150      & 6&.98e-2 & 6&.98e-2 & 13&.3 & ms & 0&.7 & ms \\
 & $h$ = 0.2, $M$ =& 10&\,000  & 4&.40e-6 & 4&.40e-6 & 517& & ms & 37&.9 & ms \\
 & $h$ = 0.2, $M$ =& 100&\,000 & 3&.27e-8 & 3&.27e-8 & 4&.8 & s & 0&.4 & s \\ \hline

\multirow{3}{*}{REXII} 
 & $h$ = 1.0, $M$ =& 38& & 2&.78e-12 & 2&.79e-12 & 10&.4 & ms & 0&.9 & ms\\
 & $h$ = 0.5, $M$ =& 65& & 1&.91e-14 & 1&.66e-14 & 15&.6 & ms & 1&.2 & ms\\
 & $h$ = 0.1, $M$ =& 278& & 7&.70e-14 & 8&.01e-14 & 51&.9 & ms & 3&.6 & ms\\ \hline
\multirow{3}{*}{RK4} 
 & \multicolumn{3}{c|}{TS = \textcolor{white}{00\,}200}    & 4&.81e-5 & 4&.87e-5   & 23&.7 & ms & 4&.8 & ms\\
 & \multicolumn{3}{c|}{TS = \textcolor{white}{0}1\,000}  & 7&.18e-8 & 7&.25e-8   & 136& & ms & 21&.4 & ms\\
 & \multicolumn{3}{c|}{TS = 50\,000} & 2&.95e-14 & 1&.63e-14 & 4&.5 & s & 1&.0 & s\\ 
\end{tabular}
\end{center}
\caption{\label{tab:w11}Comparison of accuracy and execution time between the three methods on the CPU and the GPU for wave scenario 1 \eqref{eq:WAVESCENARIO1} with a short final time $\tau=1$.}
\end{table}

From the results reported in Table~\ref{tab:w11}, we recognize that REXII outperforms the original REXI scheme and RK4 in both accuracy and execution time by a large margin. REXI is more comparable to the Runge--Kutta time stepping method of order 4 (RK4) for this initial conditions and RK4 can even outperform the original REXI method in some cases.  

For longer integration times (see the results in Table~\ref{tab:w150}), both REXI schemes drastically outperform the explicit RK4 method (as we would expect). In addition, we can see that REXII is much more accurate than the original REXI scheme.

Moreover, we observe that for all numerical methods the GPU implementation significantly outperforms the CPU implementation. For REXI and REXII the speedup ranges from approximately a factor of 7 to a factor of 15.

\begin{table}
\begin{center}
\begin{tabular}{ c| r@{\hskip 4pt}r | r@{\hskip 1pt}l | r@{\hskip 1pt}l |r@{\hskip 1pt} l@{\hskip 2pt}l |r@{\hskip 1pt} l@{\hskip 2pt}l }
\multicolumn{12}{ c }{Final time $\tau$ = 50} \\
\hline
\multirow{2}{*}{Method} & \multicolumn{2}{c}{\multirow{2}{*}{$h,M$ \hfill / \hfill time steps TS}} & \multicolumn{4}{c}{Error} & \multicolumn{6}{c}{Time}\\
& & & \multicolumn{2}{c|}{CPU} & \multicolumn{2}{c|}{GPU} & \multicolumn{3}{c|}{CPU} & \multicolumn{3}{c}{GPU} \\ \hline
\multirow{3}{*}{REXI} 
 & $h$ = 0.2, $M$ =& 7000      & 2&.63e-3 & 2&.63e-3  & 0&.44 & s & 25&  &ms \\
 & $h$ = 0.2, $M$ =& 20\,000   & 5&.11e-5 & 5&.91e-5  & 0&.98 & s & 77&  &ms \\
 & $h$ = 0.2, $M$ =& 500\,000  & 6&.35e-9 & 6&.35e-9  & 24&.1 & s & 2&.0  &s \\ \hline
\multirow{3}{*}{REXII} 
 & $h$ = 1.0, $M$ =& 1\,344 & 3&.61e-12 & 3&.41e-12 & 0&.23 & s & 18&  &ms\\
 & $h$ = 0.5, $M$ =& 2\,677 & 1&.07e-13 & 8&.93e-14 & 0&.44 & s & 46&  &ms\\
 & $h$ = 0.1, $M$ =& 13\,341& 1&.81e-13 & 1&.97e-13 & 1&.61 & s & 225&  &ms\\ \hline
\multirow{3}{*}{RK4} 
 & \multicolumn{2}{c|}{TS = \textcolor{white}{0}20\,000}   & 1&.81e-4& 1&.80e-4  & 2&.0 & s & 0&.4  &s\\
 & \multicolumn{2}{c|}{TS = 100\,000} & 2&.86e-7& 2&.84e-7 & 9&.5 & s & 2&.1  &s\\
 & \multicolumn{2}{c|}{TS = 500\,000} & 4&.56e-10& 4&.52e-10 & 47&.7 & s & 10&.3  &s\\ 
\end{tabular}
\end{center}
\caption{\label{tab:w150}Comparison of accuracy and execution time between the three methods on the CPU and the GPU for wave scenario 1 \eqref{eq:WAVESCENARIO1} with relatively large final time $\tau=50$.}
\end{table}

In the following examples we use $M$ calculated with $\rho(A)$ given by \eqref{eq:lswRoh} even if it is possible to choose a smaller one as in this example. The reason why we are doing this is that in general we can not expect the initial conditions to be that convenient. 

\subsubsection*{Wave scenario 2}

The same type of initial conditions is used as before with the exception that the frequencies are now much larger. We use:
\begin{equation}
\label{eq:WAVESCENARIO2}
\begin{aligned}
\eta(0,x,y) &= \sin(32\pi x)\cos(16\pi y) -\tfrac{1}{5}\cos(32\pi x)\sin(32 \pi y) \\
u(0,x,y) &= \cos(64\pi x)\cos(16\pi y) \\
v(0,x,y) &= \cos(32\pi x)\cos(32\pi y)
\end{aligned}
\end{equation}

\begin{table}
\begin{center}
\begin{tabular}{ c| r@{\hskip 4pt}r | r@{\hskip 1pt}l | r@{\hskip 1pt}l |r@{\hskip 1pt} l@{\hskip 2pt}l |r@{\hskip 1pt} l@{\hskip 2pt}l }
\multicolumn{12}{ c }{Final time $\tau$ = 50} \\
\hline
\multirow{2}{*}{Method} & \multicolumn{2}{c}{\multirow{2}{*}{$h,M$ \hfill / \hfill time steps TS}} & \multicolumn{4}{c}{Error} & \multicolumn{6}{c}{Time}\\
& & & \multicolumn{2}{c|}{CPU} & \multicolumn{2}{c|}{GPU} & \multicolumn{3}{c|}{CPU} & \multicolumn{3}{c}{GPU} \\ \hline
\multirow{3}{*}{REXI} 
 & $h$ = 0.2, $M$ =& 75\,003   & 2&.44e-5 & 2&.44e-5  & 3&.6 & s & 0&.27 & s \\
 & $h$ = 0.2, $M$ =& 150\,007  & 6&.35e-6 & 6&.35e-6  & 7&.1 & s & 0&.55 & s \\
 & $h$ = 0.2, $M$ =& 5\,000\,171 & 4&.52e-9 & 4&.52e-9 & 245& & s & 19&.35 & s \\ \hline
\multirow{3}{*}{REXII} 
 & $h$ = 1.0, $M$ =& 28\,448 & 4&.04e-12 & 4&.04e-12 & 3&.2 & s & 0&.47 & s\\
 & $h$ = 0.5, $M$ =& 56\,885 & 6&.53e-13 & 7&.74e-13 & 6&.3 & s & 0&.91 & s\\
 & $h$ = 0.1, $M$ =& 284\,371& 9&.36e-13 & 9&.19e-13 & 31&.6 & s & 4&.13 & s\\ \hline
\multirow{3}{*}{RK4} 
 & \multicolumn{2}{c|}{TS = \textcolor{white}{0\,}200\,000}  & 6&.11e-4 & 6&.11e-4 & 19&.7 & s & 4&.1 & s\\
 & \multicolumn{2}{c|}{TS = \textcolor{white}{0\,}500\,000}  & 1&.56e-5 & 1&.56e-5 & 49&.7 & s & 10&.1 & s\\
 & \multicolumn{2}{c|}{TS = 1\,000\,000} & 9&.77e-7 & 9&.77e-7 & 95&.8 & s & 20&.7 & s \\ 
\end{tabular}
\end{center}
\caption{\label{tab:w850} Comparison of accuracy and execution time between the three methods on the CPU and the GPU of wave scenario 2 \eqref{eq:WAVESCENARIO2} with a relatively large final time $\tau=50$.}
\end{table}
Also in this case, see Table~\ref{tab:w850}, REXII outperforms the other two methods by a large margin. Here both REXII and REXI work much better than RK4. The reason is that the high frequencies force the explicit time stepping method to take extremely small step sizes. In Table~\ref{tab:W8_M} we show that the onset of convergence strongly depends on the parameter $M$ for this problem. This is expected as $\rexii$ in the matrix case works similarly as in the scalar case. Therefore, the outcome can be compared to the results obtained in Figure~\ref{fig:exLB}.
 
\begin{table}
\begin{center}
\begin{tabular}{ c|c|c }

\multicolumn{3}{c}{Final time $\tau$ = 50, $h$ = 0.5, Method = REXII} \\
\hline
 Parameter $M$ & Error   & Time \\ 
\hline 
 $M$ = 20\,400 &  0.97     &  0.30 s\\
 $M$ = 20\,800 &  7.74e-13 &  0.32 s\\ 
\end{tabular}
\end{center}
\caption{\label{tab:W8_M}  We can observe how sensitive the choice for parameter $M$ is. For $M = 20\,400$ REXII does not approximate at all the solution, and by choosing $M$ slightly larger we obtain a precision of 12 digits. The reason why in this region convergence takes place is the same as in wave scenario 1: here the Fourier coefficients are zero for $|m|,|k| > 23$, thus in \eqref{eq:lswRoh} $D$ can be fixed to 46. }
\end{table}

\subsubsection*{Gaussian scenario}

The following initial conditions are used for the Gaussian scenario:
\begin{equation}
\label{eq:GAUSSIANSCENARIO}
\begin{aligned}
\eta(0,x,y) &= \exp(-100((x-0.5)^2 + (y-0.5)^2)) \\
u(0,x,y) &= 10^{-1}\sin(64\pi x)\sin(16\pi y)\\
v(0,x,y) &= 10^{-1}\sin(32\pi x)\sin(32\pi y)
\end{aligned}
\end{equation}

This initial function $\eta$, in contrast to the initial functions in the wave scenarios, is not exactly representable in Fourier space. Therefore, in this case we do not have the advantage of a small spectral radius or that a large part of the frequencies are equal to zero. 

The numerical results for final times $\tau=1$ and $\tau=50$ are shown in Tables~\ref{tab:G1001} and~\ref{tab:G10050}, respectively. As before, REXII outperforms the original REXI scheme and RK4 significantly in accuracy. In addition, for both REXII and the original REXI scheme the GPU implementation outperforms the CPU implementation by a factor between $7$ and $13$.

\begin{table}
\begin{center}
\begin{tabular}{ c| r@{\hskip 4pt}r@{\hskip 1pt} r | r@{\hskip 1pt}l | r@{\hskip 1pt}l |r@{\hskip 1pt} l@{\hskip 2pt}l |r@{\hskip 1pt} l@{\hskip 2pt}l }
\multicolumn{13}{ c }{Final time $\tau$ = 1} \\
\hline
\multirow{2}{*}{Method} & \multicolumn{3}{c}{\multirow{2}{*}{$h,M$ \hfill / \hfill time steps TS}} & \multicolumn{4}{c}{Error} & \multicolumn{6}{c}{Time}\\
& & & & \multicolumn{2}{c|}{CPU} & \multicolumn{2}{c|}{GPU} & \multicolumn{3}{c|}{CPU} & \multicolumn{3}{c}{GPU} \\ \hline
\multirow{3}{*}{REXI }  
 & $h$ = 0.2, $M$ =& 1&\,500    & 3&.78e-4 & 3&.78e-4  & 116& & ms & 5&.2 & ms \\
 & $h$ = 0.2, $M$ =& 3&\,000    & 3&.21e-6 & 3&.21e-6  & 224& & ms & 10&.6 & ms \\
 & $h$ = 0.2, $M$ =& 1\,000&\,025 & 4&.76e-10 & 4&.76e-10 & 48& & s & 4&.0 & s \\ \hline
\multirow{3}{*}{REXII} 
 & $h$ = 1.0, $M$ =& 580 &    & 6&.17e-13 & 6&.18e-13 & 101& & ms & 7&.2 & ms \\
 & $h$ = 0.5, $M$ =& 1\,149& & 4&.36e-15 & 6&.11e-15 & 143& & ms & 15&.2 & ms \\
 & $h$ = 0.1, $M$ =& 5\,698& & 1&.53e-14 & 1&.58e-14 & 798& & ms & 90&.6 & ms \\ \hline
\multirow{3}{*}{RK4} 
 & \multicolumn{3}{c|}{TS = \textcolor{white}{00\,}200}   & 3&.17e-2 & 3&.17e-2  & 26& & ms  & 4&.7 & ms\\
 & \multicolumn{3}{c|}{TS = \textcolor{white}{0}1\,000}  & 3&.24e-4 & 3&.24e-4  & 125& & ms  & 21&.8 & ms\\
 & \multicolumn{3}{c|}{TS = 10\,000} & 3&.13e-8 & 3&.13e-8 & 962& & ms  & 205&.0 & ms\\ 
\end{tabular}
\end{center}
\caption{\label{tab:G1001} Comparison of accuracy and execution time of the three methods on the CPU and the GPU for the Gaussian scenario \eqref{eq:GAUSSIANSCENARIO} with a short final time. }
\end{table}

\begin{table}
\begin{center}
\begin{tabular}{ c| r@{\hskip 4pt}r@{\hskip 1pt} r | r@{\hskip 1pt}l | r@{\hskip 1pt}l |r@{\hskip 1pt} l@{\hskip 2pt}l |r@{\hskip 1pt} l@{\hskip 2pt}l }
\multicolumn{13}{ c }{Final time $\tau$ = 50} \\
\hline
\multirow{2}{*}{Method} & \multicolumn{3}{c}{\multirow{2}{*}{$h,M$ \hfill / \hfill time steps TS}} & \multicolumn{4}{c}{Error} & \multicolumn{6}{c}{Time}\\
& & & & \multicolumn{2}{c|}{CPU} & \multicolumn{2}{c|}{GPU} & \multicolumn{3}{c|}{CPU} & \multicolumn{3}{c}{GPU} \\ \hline
\multirow{3}{*}{REXI} 
 & $h$ = 0.2, $M$ =& 75&\,003   & 2&.99e-6  & 2&.99e-6  & 3&.6 & s & 0&.27 & s \\
 & $h$ = 0.2, $M$ =& 150&\,007   & 8&.04e-7 & 8&.04e-7  & 7&.2 & s & 0&.55 & s \\
 & $h$ = 0.2, $M$ =& 5\,000&\,171  & 5&.77e-10 & 5&.77e-10 & 239& & s & 19&.35 & s \\ \hline

\multirow{3}{*}{REXII} 
 & $h$ = 1.0, $M$ =& 28&\,448 & 6&.18e-13 & 6&.46e-13 & 3&.4 & s & 0&.47 & s\\
 & $h$ = 0.5, $M$ =& 56&\,885 & 6&.06e-14 & 6&.90e-14 & 6&.5 & s & 0&.91 & s\\
 & $h$ = 0.1, $M$ =& 284&\,371& 1&.04e-13 & 4&.48e-14 & 31&.0 & s & 4&.13 & s\\ \hline
\multirow{3}{*}{RK4} 
 & \multicolumn{3}{c|}{TS = \textcolor{white}{0\,}200\,000} & 6&.06e-5 & 6&.06e-5  & 18&.4 & s  & 4&.14 & s\\
 & \multicolumn{3}{c|}{TS = \textcolor{white}{0\,}500\,000} & 1&.54e-6 & 1&.54e-6  & 47&.6 & s & 10&.16 & s\\
 & \multicolumn{3}{c|}{TS = 1\,000\,000} & 9&.67e-8 & 9&.67e-8 & 95&.6 & s & 20&.57 & s\\ 
\end{tabular}
\end{center}
\caption{\label{tab:G10050} Comparison of accuracy and execution time of the three methods on the CPU and the GPU for the Gaussian scenario \eqref{eq:GAUSSIANSCENARIO} with a long final time.}
\end{table}

\section{Conclusion}

The original REXI scheme is already a good method to compute the action of the matrix exponential parallel in time. The main downside is that it is not very precise. We proposed a modification of the REXI approach that achieves accuracy close to machine precision at similar or, for some problems, even lower computational cost. The strength of the $\rexi$ methods is the fact that they can be easily parallelized in time (in addition to the commonly used parallelization in space). We have demonstrated this by providing an implementation on massively parallel graphic processing units. The GPU implementation shows a drastic speedup compared to the CPU implementation. 

\section{Acknowledgements}

This project has received funding from the European Union’s Horizon 2020 research and innovation programme under the Marie Skłodowska-Curie grant agreement No 847476. The views and opinions expressed herein do not necessarily reflect those of the
European Commission.

\bibliographystyle{plain}
\bibliography{literatur}

\section*{Appendix A}

\begin{table}[H]
\begin{small}
\begin{center}
  \begin{tabular}{ l l@{\hskip 1pt} l}
  \hline
     $\mu=$&-&\texttt{5.133333333333333}\\
    \hline
	$a_{0}=$&-&\texttt{6.520430828919864e+01}                              \\
	$a_{1}=$& &\texttt{4.261818064131437e+01 + 2.761406741120911e+01i}             \\
	$a_{2}=$&-&\texttt{9.801650304425239e+00 - 2.189295463610722e+01i}           \\
	$a_{3}=$&-&\texttt{1.054225194693395e+00 + 6.791786454153551e+00i}            \\
	$a_{4}=$& &\texttt{7.950505668209775e-01 - 8.904997258367445e-01i}         \\
	$a_{5}=$&-&\texttt{1.218558380859130e-01 + 3.321241563407446e-02i}          \\
	$a_{6}=$& &\texttt{7.365401806949337e-03 + 2.212802103193251e-03i}       \\
	$a_{7}=$&-&\texttt{2.801087265991056e-04 - 5.566945197754387e-04i}   \\
	$a_{8}=$& &\texttt{1.254835436432561e-04 - 2.467200513365371e-04i}    \\
	$a_{9}=$& &\texttt{2.295472292491263e-04 - 8.494118951459107e-05i}    \\
	$a_{10}=$& &\texttt{1.858484460459430e-04 + 9.242889460185034e-05i}     \\
	$a_{11}=$& &\texttt{4.068056518449676e-05 + 1.653479957565515e-04i}    \\
	$a_{12}=$&-&\texttt{8.341508001647741e-05 + 1.045331460447588e-04i}   \\
	$a_{13}=$&-&\texttt{9.970528169841103e-05 - 5.856228484297677e-06i}  \\
	$a_{14}=$&-&\texttt{3.499639858693093e-05 - 6.129059473910835e-05i} \\
	$a_{15}=$& &\texttt{2.295021920298455e-05 - 4.099832469456381e-05i}   \\
	$a_{16}=$& &\texttt{2.931048772724314e-05 + 1.708815129697846e-07i}   \\
	$a_{17}=$& &\texttt{7.502088478301169e-06 + 1.525082051744077e-05i}    \\
	$a_{18}=$&-&\texttt{5.815291167450100e-06 + 6.919604247338349e-06i}     \\
	$a_{19}=$&-&\texttt{4.069948458364005e-06 - 1.440010113050771e-06i}  \\
	$a_{20}=$& &\texttt{7.932524475429588e-08 - 1.794169428574330e-06i}    \\
	$a_{21}=$& &\texttt{6.120984882186265e-07 - 1.131894636585849e-07i}   \\
	$a_{22}=$& &\texttt{5.531365159161319e-08 + 1.585749903175946e-07i}    \\
	$a_{23}=$&-&\texttt{2.867805871375946e-08 + 1.239499740327838e-08i}   \\
	$a_{24}=$&-&\texttt{1.143081277095316e-09 - 2.763239274253499e-09i}  \\
    \hline
  \end{tabular}
\end{center}
\end{small}
\caption{\label{tab:coef_al} Coefficients $a_l$ for \eqref{eq:ratapproxgaus} with $a_l = \overline{a_{-l}}$.}
\end{table}

\section*{Appendix B}

In \cite{HBMW2015} the authors give the idea how to determine an upper error bound, which we will analyze here more in depth. The following expression has to be estimated:

\begin{equation}
\label{eq:REXIerror}
\Biggl|e^{ix} - \sum_{m = -M}^M b_m \sum_{l = -L}^L \re\left(\frac{ha_l}{ix+h(\mu + i(m+l))}\right)\Biggr| .
\end{equation}
To do so, the following sum is added and subtracted inside the modulus
\[
\sum_{m = -M}^M b_m \psi_h(x+mh)
\]
and then the triangle inequality is applied. Thus, we end up with the two terms
\begin{equation}
\label{eq:REXIerror1}
\Biggl|e^{ix} - \sum_{m = -M}^M b_m \psi_h(x+mh)\Biggr| 
\end{equation}
and
\begin{equation}
\label{eq:REXIerror2}
\Biggl| \sum_{m = -M}^M b_m \sum_{l = -L}^L \re\left(\frac{ha_l}{ix+h(\mu + i(m+l))}\right) - \sum_{m = -M}^M b_m \psi_h(x+mh)\Biggr|.
\end{equation}
The term \eqref{eq:REXIerror2} can be estimated using \eqref{eq:bm}:
\begin{align*}
\Biggl| \sum_{m = -M}^M b_m&  \left(\sum_{l = -L}^L \re\left(\frac{ha_l}{ix+h(\mu + i(m+l))}\right) - \psi_h(x+mh) \right)\Biggr| \\
&\leq \sum_{m = -M}^M e^{h^2}  \Biggl|  R\left(\frac{x}{h}+m\right) - \psi_h(x+mh) \Biggr| \\ 
&\leq e^{h^2}(2M+1)\delta_{2},
\end{align*}
where $\delta_{2}$ is the approximation error of the Gaussian function. In our case, $\delta_2 = 8\times 10^{-15}$, see step 2 in Section~\ref{sec:REXI_derivation_scalar}. 
This bound for \eqref{eq:REXIerror2} is not sharp. It depends linearly on the parameter $M$, but we observed that after doing several numerical experiments this amount does not increase much when $M$ increases.

To deduce an error bound for \eqref{eq:REXIerror1}, Poisson's summation formula is applied: 
\[
\sum_{m = -\infty}^\infty \phi_h(x + mh) = 
 \frac{1}{h}\sum_{k = -\infty}^\infty e^{2\pi i(k/h)x} \hat{\phi}_h \left(\frac{k}{h}\right).
\]
For the choice  
\[
\phi_h(x+mh) = e^{-i(x+mh)}\psi_h(x+mh)
\]
we get
\begin{align*}
e^{-ix}\sum_{m = -\infty}^\infty e^{-imh}\psi_h(x+mh) 
 = \frac{1}{h}\sum_{k = -\infty}^\infty e^{2\pi i(k/h)x} \hat{\psi}_h \left(\frac{k}{h} +\frac{1}{2\pi}\right).
\end{align*}
This equation is then multiplied by $e^{ix}$ and the sum on the right-hand side of the equation is split up, leading to

\begin{align*}
\sum_{m = -\infty}^\infty e^{-imh}\psi_h(x+mh) -e^{ix}\frac{1}{h}\hat{\psi}_h\left(\frac{1}{2\pi}\right) 
 &= \frac{1}{h}\sum_{k \neq 0} e^{2\pi ix(\frac{k}{h} + \frac{1}{2\pi} )} \hat{\psi}_h \left(\frac{k}{h} +\frac{1}{2\pi}\right). \\
\end{align*}
This implies

\begin{align*}
\Biggl|\sum_{m = -M}^M &e^{-imh} \psi_h(x+mh) -e^{ix}\frac{\hat{\psi}_h\left(\frac{1}{2\pi}\right)}{h}\Biggr| = \\
 &= \Biggl| \frac{1}{h} \sum_{k \neq 0} e^{2\pi ix(\frac{k}{h} + \frac{1}{2\pi} )} \hat{\psi}_h \left(\frac{k}{h} +\frac{1}{2\pi}\right) 
-\sum_{|m|>M} e^{-imh} \psi_h(x+mh) \Biggr| \\
 &\leq \frac{1}{h}\sum_{k \neq 0} \hat{\psi}_h \left(\frac{k}{h} +\frac{1}{2\pi}\right) + \sum_{|m|>M} \psi_h(x+mh). \\
\end{align*}
The coefficients $b_m = he^{-imh}\hat{\psi}_h\left(\frac{1}{2\pi}\right)^{-1}$ are given by~\eqref{eq:bm} and $\hat{\psi}_h(\omega)$ by~\eqref{eq:psi_fourier}. Thus, if  the inequality is divided by the positive number $\frac{\hat{\psi}_h\left(\frac{1}{2\pi}\right)}{h}$ it follows

\begin{equation}
\label{eq:REXIerrdelta1} 
\begin{aligned}
\Biggl|\sum_{m = -M}^M& b_m \psi_h(x+mh) -e^{ix}\Biggr| \leq \\
 &\leq \frac{1}{\hat{\psi}_h\left(\frac{1}{2\pi}\right)}\left(\sum_{k \neq 0} \hat{\psi}_h \left(\frac{k}{h} +\frac{1}{2\pi}\right) + h \sum_{|m|>M} \psi_h(x+mh)\right) \\
&\leq e^{h^2}\left(\frac{1}{h}\sum_{k \neq 0} \hat{\psi}_h \left(\frac{k}{h} \right) +  \sum_{|m|>M} \psi_h(x+mh)\right) \\
&= e^{h^2}\left(\sum_{k \neq 0} e^{-4\pi^2k^2}+\sum_{|m|>M} \psi_h(x+mh)\right):= e^{h^2}\delta_1.
\end{aligned} 
\end{equation}
So the final result is the following:
\begin{equation}
\label{eq:REXIerrStim}
|e^{ix} - \rexi(ix,M,h)| \leq e^{h^2}(\delta_1 + (2M+1)\delta_2).
\end{equation}
Let us now analyse $\delta_1 = \sum_{k \neq 0} e^{-4\pi^2k^2}+\sum_{|m|>M} \psi_h(x+mh)$. The first sum is neglible. Since $e^{-4\pi^2k^2} \leq e^{-4\pi^2k}$ it follows that 
\[
\sum_{k = 1}^{\infty} e^{-4\pi^2k^2} < \sum_{k = 1}^{\infty} e^{-4\pi^2k} = \frac{1}{1-e^{-4\pi^2}}-1 \approx 7.15 \times 10^{-18}.
\]

The second sum is more interesting regarding the error. It can be shown that for $|x| \leq (M-m_0)h$ this sum is also neglible, where $m_0$ is a constant related to how small the first two terms of the sum $\sum_{|m|>M}\psi_h(x+mh)$ are. An easy calculation shows that 

\[
\psi_h(z) \leq \text{tol} \quad \text{if} \quad |z| \geq 2h\sqrt{-\log(\sqrt{4\pi}\text{ tol})} = ch.
\]
With $\text{tol} = 10^{-16}$ we have $c \approx 12$. Now it is of interest that $\psi_h(x + (M+1)h) < \text{tol}$ and $\psi_h(x -(M+1)h) < \text{tol}$ such that these two shifted Gaussian functions have no impact on the approximation of $e^{ix}$. This implies, by the previous calculation, that $(x - (M+1)h) \leq -ch$ and $(x +(M+1)h) \geq ch$. Thus, $|x| \leq (M+1-c)h$ and we may define $m_0 := c-1 \approx 11$. 

\end{document}